\documentclass[a4paper]{article}

\usepackage[latin1]{inputenc} %
\usepackage[T1]{fontenc} %
\usepackage{RR}
\usepackage{hyperref}

\usepackage{graphicx,amsmath,amssymb}
\usepackage{subfigure}
\usepackage{epsfig}
\usepackage{verbatim}
\usepackage{algorithm}
\usepackage{algorithmic}
\usepackage{color}
\usepackage{multirow}
\usepackage{latexsym}

\newcommand{\IR}{\mathbb{R}}
\renewcommand{\d}{\mathrm{d}}
\newcommand{\D}{\mathrm{D}}
\newcommand{\p}{p}
\newcommand{\n}{n}
\renewcommand{\pm}{p_m}
\newcommand{\parder}[2]{\frac{\partial #1}{\partial #2}}
\newcommand{\e}{\mathrm{e}}
\newcommand{\ch}{\mathrm{ch}}

\newcommand{\sU}{\mathsf{U}}
\newcommand{\cT}{\mathcal{T}}
\newcommand{\tV}{U}

\newcommand{\Om}{\Omega}
\newcommand{\om}{\omega}

\RRdate{June 2010}

\RRauthor{%
P. Bernhard\thanks{INRIA Sophia Antipolis-M\'editerran\'ee, France, Pierre.Bernhard@sophia.inria.fr}
\and
F. Grognard\thanks{INRIA Sophia Antipolis-M\'editerran\'ee, France, Frederic.Grognard@sophia.inria.fr}
\and
L. Mailleret\thanks{INRA Sophia Antipolis, France, Ludovic.Mailleret@sophia.inra.fr}
\and
A.R. Akhmetzhanov\thanks{INRIA Sophia Antipolis-M\'editerran\'ee, France, akhmetzhanov@gmail.com}
}

\authorhead{Bernhard \& Grognard \& Mailleret \& Akhmetzhanov}

\RRtitle{Stabilité évolutionnaire des stratégies d'histoires de vies de consommateurs solidaires confrontés  à une invasion de mutants}
\RRetitle{ESS for life-history traits of cooperating consumers facing cheating mutants}
\titlehead{ESS for life-history traits of cooperating consumers facing cheating mutants}

\RRresume{
Nous considérons une population d'individus identiques s'attaquant à une ressource épuisable. Les individus de la population doivent choisir une stratégie qui définit comment ils utilisent leur temps disponible au cours de leur vie à se nourrir, pour la
reproduction (par exemple la ponte), ou en partageant leur énergie entre ces deux activités.
Nous supposons ici que leur vie ne dure qu'une saison complète, de sorte que la stratégie choisie se traduit par la production d'un certain nombre d'oeufs pondus au cours de la
saison. Ce nombre permet alors de définir l'évolution à long terme de la population car les oeufs constituent la base de la population pour la saison suivante.
Dans ce rapport, nous nous concentrons strictement sur ce qui se produit au sein d'une saison en considérant deux stratégies possibles: l'optimum collectif et la stratégie non invasible.

La stratégie (collective) optimale définit un arc singulier le long duquel la ressource est économisée; sur le long terme, cette stratégie peut conduire à un équilibre. Par contre, cettre stratégie rend la population sensible à l'invasion par un mutant glouton qui n'économiserait pas la ressource sur l'arc singulier. Par conséquent, la stratégie ``optimale'' n'est pas évo\-lu\-tion\-nai\-re\-ment stable.

Nous allons donc rechercher une stratégie évolutionnairement stable, en utilisant la première condition d'ESS (``condition de Nash-Wardrop''). Il en résulte un autre arc singulier
et une stratégie qui préserve moins la ressource.

Dans cet article, nous vérifions aussi si une stratégie mixte sur un arc singulier peut être interprétée comme une représentation d'une population dont une partie se reproduit continuellement et l'autre se nourrit. La réponse est négative: la stratégie mixte ne peut représenter des parts de populations jouant des stratégies pures.
}

\RRabstract{
We consider a population of identical individuals preying on an exhaustible
resource. The individuals in the population choose a strategy that defines
how they use their available time over the course of their life for feeding, for
reproducing (say laying eggs), or split their energy in between these two activities.
We here suppose that their life lasts a full season, so that the chosen
strategy results in the production of a certain number of eggs laid over the
season. This number then helps define the long-run evolution of the population
since the eggs constitute the basis for the population for the next season.
However, in this paper, we strictly concentrate on what occurs within a season,
by considering two possible strategies: the collective optimum and the
uninvadable strategy.

The (collective) optimal strategy involves a singular arc, which has some
resource saving feature and may, on the long-term, lead to an equilibrium.
However, it is susceptible to invasion by a greedy mutant which free-rides
on this resource saving strategy. Therefore, the ``optimal'' strategy is not
evolutionarily stable.

We thus look for an evolutionarily stable strategy, using the first ESS
condition (``Nash-Wardrop condition''). This yields a different singular arc
and a strategy that conserves less the resource. Unfortunately, we show in a
forthcoming paper that, though it is an ESS, it is not long-term (interseasonal)
sustainable. This is an instance of the classical ``tragedy of the commons''.
In this paper we also investigate whether the mixed strategy on the singular
arcs may be interpreted in terms of population shares using pure strategies.
The answer is negative.
}

\RRmotcle{Jeux différentiels, Dynamique des populations, Modèle consommateur-ressource, ESS}
\RRkeyword{Differential Games, Population Dynamics, Resource-Consumer Model, Evolutionarily Stable Strategy}

\RRprojet{Comore}  
\RRtheme{\THBio} 
\URSophia 

\begin{document}

\RRNo{7314}

\makeRR

\section{Introduction}
\label{Section_Introduction}
\subsection{Context and contribution}
\subsection{Context and contribution}

The evolution of life history traits, {\it e.g.} the reproduction
strategy an organism should follow during the course of its life, is
an important topic of evolutionary ecology
\cite{Lessells1991,Stearns1992}. In this respect, seminal
contributions based on optimal control theory principles have focused
on the optimal energy allocation schedules of organisms into growth or
reproduction (see the reviews
\cite{Perrin1993,Heino1999,Iwasa2000}). These studies consider an
``energy allocation'' trait, the proportion of energy an organism
invests at a given Â moment of its life into somatic growth or into
reproduction forms, and look for dynamic (age and state dependent)
strategies that maximize the reproductive output of the organism. They
show in particular how different growth-reproduction dynamical
strategies, like the determinate and indeterminate growth patterns,
may emerge from different hypotheses on the biology of the organisms
(\cite{Perrin1993,Heino1999}). These studies are however based on
models at the individual level, and thus do not address population
issues and the potential feedbacks of population densities on the life
history strategies that should be promoted by evolution.

The eventual success of a strategy does however not necessarily depend on its
optimal exploitation of their ressource: a species indeed participates in games with other species and mutants of its own which could take advantage of its ``optimal'' strategy; evolutionary game theory is then the general framework for the study of which strategy is selected through evolution. In that framework, the Adaptive Dynamics approach Â is a recent methodology specifically
tailored to tackle evolutionary questions arising in a population
context \cite{Geritz1997,Dercole2008}. Based on the (un)-invasibility
analysis of rare mutants into a resident population at ecological
equilibrium, the studies adopting this framework have mostly
concentrated on scalar or multi-dimensional traits. As such, they
fairly overlooked potential dynamical life history traits which are by
essence of an infinite-dimensional nature (see however
\cite{Eskola2009} for an exception).

In this contribution, we shall build on the model proposed by
\cite{Akhmetzhanov2010} to describe the Â interactions between
populations of seasonal species of consumers and resource
organisms. This model is of a semi-dicrete type \cite{Mailleret2009}:
the intra-seasonal interactions between the mature consumers and
resource populations as well as production of immature forms are
modelled with a continuous time system, while the inter-seasonal
demographic processes, essentially maturation of the immature forms
and death of the mature populations, are modelled with a discrete
mapping. As such, this model allows to account for dynamical
strategies ({\it e.g.} feeding or reproduction of the consumers) at
the intra-seasonal level, as well as their consequences on populations
at the long term scale through the investigation of the inter-seasonal
dynamics. More precisely, \cite{Akhmetzhanov2010} studied the optimal
feeding-reproduction strategy of the consumer population; it can be
looked at as an extension of the works on optimal energy allocation
described above to an explicit populational context. The current
consensus is however that evolutionary studies should rely on
invasibility criteria (as the Adaptive Dynamics methodology does)
rather than on optimization criteria \cite{Mylius1995, Metz2008}.

The objective of the present contribution is to extend the work of
\cite{Akhmetzhanov2010} to investigate ``un-invadable'' dynamic
strategies, {\it i.e.} strategies that can not be beaten by
alternative strategies followed by some mutants taking advantage of
their small density. In particular, we show that the optimal
cooperating strategy of the consumers computed in
\cite{Akhmetzhanov2010} can be invaded (see also
\cite{Akhmetzhanov2010a} for a deeper study of this issue), and that
consumers may adopt an alternative Evolutionarily Stable Strategy
(ESS) which prevents any mutant invasion. Though it will not be developed
in the present paper which concentrates on within season strategies, we can
show that, contrarily to the ``cooperative'' or ``optimal'' strategy that can
lead to a long term co-existence of resource and consumer populations,
the ESS strategy always leads to an
evolutionary suicide of the consumers in the interseasonal dynamics.

\subsection{The model}
We consider an homogeneous population of essentially fixed size.
The energy level of a representative individual is $\p \in \IR_+$, the
proportion of its time spent feeding is $u\in[0,1]$, the amount of
available resource (say the preys population size) is $\n \in \IR^+$.
The dynamics are as follows, with $a$, $b$ and $c$ positive coefficients:
\begin{eqnarray}
\dot\p &= &-a\p + b\n u\,, \label{pdot}\\
\dot\n &= &-c\n u\,. \label{ndot}
\end{eqnarray}
It should be noticed that these dynamics leave the strictly positive othant
invariant. We may therefore exploit the following feature of these dynamics:
The ratio $x = \p/\n$ obeys the scalar dynamics
\begin{equation}\label{xdot}
\dot x = -ax + (b+cx)u\,.
\end{equation}

The criterion for long term sustainability is the number of eggs layed during
a season, of length $T$, given as
\begin{equation}\label{J}
J = \int_0^T\!(1-u(t))\p(t)\,\d t\,.
\end{equation}
This criterion is \emph{not} expressible in terms of $x$ only.

A small population of potential invaders will be considered, with the same
energy dynamics
\begin{equation}\label{pmdot}
\dot\p_m = -a\pm + bnu_m
\end{equation}
and criterion
\begin{equation}\label{Jm}
J_m=\int_0^T\!(1-u_m)\pm\,\d t\,,
\end{equation}
but with no effect on the resource amount because of the negligible size
of that population. It gives rise to the variable $x_m=\pm/\n$ with dynamics
\begin{equation}\label{xmdot}
\dot x_m = -ax_m + bu_m+cx_m u\,.
\end{equation}

\section{Collective optimum}
\subsection{Euler Pontryagin equations}
We use a standard maximum principle approach to finding the optimal
collective behavior. To take advantage of the reduced size dynamics,
let $V$ be the Bellman function of this problem. We shall show
that it has a solution of the form
\[
V(t,\n,\p) = \n \tV (t,x)\,.
\]
If this is so,
\[
\parder{V}{t} = \n\parder{\tV }{t}\,,\quad
\parder{V}{\p} = \parder{\tV }{x}\,,\quad
\parder{V}{\n} = \tV (t,x) - x\parder{\tV }{x}\,.
\]
The Hamilton Jacobi Caratheodory Bellman equation thus reads
\[
\n\parder{\tV }{t} + \max_{u\in[0,1]}\left[-\parder{\tV }{x}(a\p-b\n u)
-\left(\tV  - x\parder{\tV }{x}\right)c\n u + (1-u)\p \right] = 0\,.
\]
Again, we may divide through by $\n$ to obtain
\begin{equation}\label{HJCB}
\parder{\tV }{t} + \max_{u\in[0,1]}\left[ -\parder{\tV }{x}[ax-(b+cx)u]
-\tV cu + (1-u)x\right] = 0\,,
\end{equation}
with the terminal condition
\[
\forall x\in\IR_+\,,\quad \tV (T,x) = 0\,.
\]

Introduce the notations
\[
\nu(t) = \parder{\tV }{t}(t,x(t))\,,\quad \lambda(t) = \parder{\tV }{x}(t,x(t))\,.
\]
The equations of the characteristics of this PDE, or the Euler-Pontryagin
adjoint equations, yield
\begin{equation}\label{EPl}
\begin{array}{ll}
\dot\nu = cu\nu\,, & \nu(T) = 0\,,\\
\dot\lambda = a\lambda -1 + u\,, & \lambda(T) = 0\,,
\end{array}
\end{equation}
with $u \in \mathrm{argmax}_u\{[\lambda(b+cx)-x-c\tV ]u\}$.
We further simplify this with the additional notation
\(
\mu = \tV  - \lambda x
\)
which is just $\partial V/\partial n$. Using (\ref{EPl}) and (\ref{HJCB}), we get
\begin{equation}\label{EPm}
\dot\mu = (c\mu-b\lambda)u\,,\quad \mu(T)=0\,,
\end{equation}
and that the optimal $u$ is given by
\begin{eqnarray}
u(t) \!&\in &\!\mathrm{argmax}_u\{\sigma(t)u\}\,, \label{uopt}\\
\sigma \!&= &\! b\lambda - c\mu - x\,. \label{sigma}
\end{eqnarray}

\subsection{Field of trajectories}
\subsubsection{The switch line}
At $t=T$, $\sigma(T) = -x(T) < 0$. Therefore, the final optimal control is $u=0$.
The final part of the optimal state and adjoint trajectories are thus
\[
x(t) = x(T)\e^{a(T-t)}\,, \quad \lambda(t)=\frac{1}{a}\left(1- \e^{-a(T-t)}\right)\,, \quad
\mu(t)=0\,,
\]
hence $\sigma(t) = b\lambda(t) - x(t)$ vanishes on a \emph{switch line}
\[
x=b\lambda=(b/a)[1-\e^{-a(T-t)}]\,.
\]
On this curve, the optimal control switches
from 1 to 0.

A terminal trajectory with $u=0$ can emanate from this switch line
only if its slope is larger than that of the switch line itself, i.e. if
\[
-ax \ge -b\e^{-a(T-t)}
\]
which, placed back in the equation of the switch line yields
\begin{equation}\label{bar}
x\le \hat x = b/(2a) \quad\mbox{and}\quad t \ge \hat t = T-(\ln 2)/a\,.
\end{equation}
Therefore, the switch line extends only from this $(\hat t, \hat x)$ to $x(T)=0$.
The last primary trajectory, passing through $x(\hat t)=\hat x$ has $x(T)=b/(4a)$.
The trajectories with $x(T) > b/(4a)$ have $u=0$ as their optimal control for
all $t\le T$.

\subsubsection{Tributaries of the switch line}
Trajectories with $x(T) < b/(4a)$ have a corner on the switch line. Prior to
reaching it, they have $u = 1$, hence
\[
\dot x = (c-a)x + b\,.
\]
These trajectories, reffered to as the tributaries of the switch line, are as follows.
Let $t_s$, $x_s=(b/a)[1-\e^{-a(T-t_s)}]$ be their end point on the switch line.
We get
\[
x(t) = \left(x_s+\frac{b}{c-a}\right)\e^{-(c-a)(t_s-t)} - \frac{b}{c-a}\,.
\]
We show in the appendix that indeed, $\sigma$ does not change sign on
these trajectories. They have $x(\cdot)$ increasing with time. As a matter of
fact, we have $\dot x > 0$ if $c\ge a$ or, if $c < a$ and $x \le b/(a-c)$.
But since on the switch line, $x \le b/(2a) < b/(a-c)$, this is always satisfied.
Hence all such trajectories come from the time axis. However, it is easy to
check that, if
\[
T > \frac{1}{a}\ln 2 + \frac{1}{c-a}\ln\frac{c+a}{2a}\,,
\]
the trajectory ending at $x(\hat t) = \hat x$ comes from $x(t)=0$ for a
positive $t$, and all other, which have a smaller $x(t)$, a fortiori. Hence
in that case none of the trajectories described so far accounts for an initial
condition $x(0) \le (b/4a)\e^{aT}$.

\subsubsection{Singular line and tributaries}
The space between the primary trajectory $x(t) = b/(4a)\exp[a(T-t)]$ and
the tributary
\[
x(t) = b\left(\frac{1}{2a} +\frac{1}{c-a}\right) \e^{-(c-a)(t_s-\hat t)} - \frac{b}{c-a}
\]
is not accounted for so far. In that space, we construct a singular line
ending at $(\hat t, \hat x)$. We differentiate twice the expression for $\sigma$
(\ref{sigma}), using (\ref{xdot}),(\ref{EPl}),(\ref{EPm}):
\[
\begin{array}{l}
\dot\sigma = c\sigma + ab\lambda + ax - b\,,\\
\ddot\sigma = c\dot\sigma +a[(2b+cx)u + ab\lambda -b - ax]\,.
\end{array}
\]
Using $\sigma = 0$ and hence $\dot\sigma=ab\lambda+ax-b=0$, we get
\[
\frac{1}{a}\ddot\sigma = (2b+cx)u - 2ax\,.
\]
Hence, a singular control is
\begin{equation}\label{usingular}
u = \frac{2ax}{2b + cx}\,.
\end{equation}
With this control,
\[
\dot x = \frac{acx^2}{b+2cx} > 0\,.
\]
Hence in particular $x(t) < \hat x = b/(2a)$. And towards the negative times,
the asymptote is $x=0$. More precisely, the singular arc satisfies
\[
\ln\frac{\hat x}{x} + \frac{2b}{c}\left(\frac{1}{\hat x}-\frac{1}{x}\right) = a(\hat t-t)\,,
\]
which can be solved for $x$ in terms of $t$ with the help of the Lambert function
$W $ as
\[
x=\frac{2b}{c}\left[ W_0\!\left(\frac{2a}{c}\e^{a(T-t)+4a/c}\right)\right]^{-1}\,,
\]
a perfectly useless formula \ldots

We check that (\ref{usingular}) yields $u \le 1$. This is true if $(2a-c)x \le 2b$.
This is always true if $2a \le c$, and, if $2a > c$ provided that $x \le b/(2a-c)$.
But we have seen that on the singular arc, $x < b/(2a) < b/(2a-c)$. Hence the
conditionn is always satisfied.

Substituting (\ref{usingular}) into (\ref{EPl}) and (\ref{EPm}), and using (\ref{bar}),
one obtains
\[
\begin{array}{ll}
\displaystyle{\dot\lambda = -\frac{1}{b}\dot x\,, }&
\displaystyle{\lambda = Â \frac{1}{a}-\frac{x}{b}\,,}\\ [9pt]
\displaystyle{\dot\mu = -\frac{2}{c}\dot x\,,} &
\displaystyle{\mu = \frac{b}{ac}-\frac{2x}{c}}\,.
\end{array}
\]

Now, tributaries can reach this singular arc at any of its points with either $u=1$,
coming from ``under'' the singular arc, or with $u=0$ coming from ``above''.
Let $(t^s,x^s)$ be a point where a tributary reaches this singular line. The
tributary with $u=0$ satisfy
\[
\begin{array}{l}
\displaystyle{x = e^{a(t^s-t)}x^s\,,}\\ [6pt]
\displaystyle{\lambda = \frac{1}{a}-\frac{x^s}{b}e^{-a(t^s-t)}}\,,\\ [12pt]
\displaystyle{\mu = \frac{1}{c}\left(\frac{b}{a}-2x^s\right)}\,,
\end{array}
\]
and hence
\[
\sigma = 2x^s(1-\ch[a(t^s-t)]) \le 0\,.
\]
Therefore, the optimal control remains $u=0$ all along. The analysis of the
tributaries with $u=1$ is done in the appendix together with those of the
switch line.

\section{Uninvadable strategy}
\subsection{``Optimal'' strategy can be invaded}
We consider a very small sub-population of mutating individuals behaving
differently. Let a subscript $m$ stand for variables concerning that population,
with $x_m=\pm/n$. It is so small that it has no effect on resource depletion.
Its dynamics are given by (\ref{pmdot}) and (\ref{xmdot}), while the dynamics
(\ref{pdot}) and (\ref{ndot}) are unchanged, and therefore also (\ref{xdot}).
Its fitness is measured by the criterion (\ref{Jm}).

Let us consider $u = \phi(t,x)$ given, and $u_m$ as a control seeking to
maximize $J_m$. The same construction as above leads to the Value function
$V_m(\pm,n)\tV _m(t,x_m)$, with the Hamilton Jacobi equation
\begin{equation}\label{HJm}
\parder{\tV _m}{t} -\parder{\tV _m}{x_m}(ax_m-bu_m+cx_m u)
-\tV _mcu + (1-u_m)x_m = 0\,.
\end{equation}
It should be pointed out that this equation holds true with $V_m = n\tV _m(t,x_m)$
the fitness obtained by the mutant, even if it does not ``optimize''.
Using the notations
\[
\nu_m = \parder{\tV _m}{t}\,, \quad \lambda_m=\parder{\tV _m}{x_m}\,,
\quad \mu_m = \tV _m-\lambda_mx_m\,,
\]
its characteristics are
\begin{equation}\label{ELm}
\begin{array}{ll}
\dot \nu_m = cu\nu_m\,, & \nu_m(T)=0\,,\\
\dot\lambda_m = a\lambda_m-1+u_m\,, & \lambda_m(T)=0\,,\\
\dot\mu_m = c\mu_mu-b\lambda_mu_m\,.
\end{array}
\end{equation}
Differentiating with respect to $u_m$, we find that
\[
\D J_m.v_m = \int_0^T\!\sigma_m(t)v_m(t)\,\d t\,.
\]
with
\begin{equation}\label{sigmam}
\sigma_m = b\lambda_m-x_m\,.
\end{equation}

If the mutant behaves as the resident population, the variables indexed $m$
are equal to the corresponding variables without the index. On the
singular arc, this leads to $\sigma_m = \sigma + c\mu$. We also have, on
this arc, using $\sigma=b\lambda-c\mu-x=0$,
\[
\dot\mu = (c\mu-b\lambda)u=-\frac{2ax^2}{2b+cx} < 0\,.
\]
Since $\mu(\hat t)=0$, it follows that for $t < \hat t$,
$\mu(t)>0$, hence $\sigma_m>0$,
which means that the mutant can improve its fitness by choosing $u_m=1$.
Hence a mutant may appear that will be more fit than the resident population:
an invasion.

This is not surprising: the singular arc corresponds to a time interval during
which the resident population cooperates to spare some resource. The mutant,
who does not deplete that resource, has no incentive to spare it and can use
as much as it wants. We defer to a future paper the analysis of the optimum
behaviour of the mutant in this cooperating ``optimum'' population.

\subsection{An uninvadable strategy}\label{uninvadable}
Since the ``optimal'' strategy can be invaded by a mutant, there is a priori no reason for this strategy to survive on the long-run. Therefore, we seek a strategy for the resident population that could not be invaded.
For that purpose, we seek a control law $u=\phi(x,t)$ such that the optimum
mutant strategy will be
\begin{equation}\label{equal}
u_m=u\,.
\end{equation}
Therefore, there will be no possibility
for a mutant to be more fit than the resident. This will lead to a new type of singular arc.

Close to final time, we have $\sigma_m=0$, therefore the optimal behaviour
of the mutant is $u_m=0$. We choose $u=0$ to insure the equality (\ref{equal}).
We obtain the same field as previously, from the switch line on, up to time $T$.

We attempt to construct a singular line for the mutant, attaching to
$(\hat t,\hat x)$, and satisfaying (\ref{equal}). This means that the dynamics
of both the state and the adjoints are as in the previous analysis, because of
(\ref{equal}), but $\sigma_m$ is as in (\ref{sigmam}). Differentiating along the
singular arc, we get
\[
\dot\sigma_m = ba\lambda -b+ax-cxu
\]
hence, $\sigma_m=0$ translates into
\begin{equation}\label{singularm}
u = \frac{ax+ab\lambda-b}{cx}= \frac{2ax-b}{cx} =: u^\sigma(x)\,.\,.
\end{equation}
The second form above is obtained using $\sigma_m=0$.
Interestingly, the first differentiation of $\sigma_m$ yields a definite expression
for $u$, contrary to what is observed in classical singular arcs. The reason is
that this $u$ is not the optimizing one, say $u_m$, but the control of the
resident population. However, we wish that both be equal. This also leads to
$x_m=x$. Yet, the formula (\ref{singularm}) should not be considered a feedback
\[
u = \frac{2ax_m-b}{cx_m}
\]
which would cause $x_m$, and therefore $u_m$, to influence the mutant's dynamics,
resulting in a different set of adjoint equations. The resident's control and state
histories should be considered as fixed a priori, independant of the mutant's action.
They are just chosen in such a way that the optimal mutant's variables coincide
with them.

Substituting $u^\sigma$ (\ref{singularm}) in the dynamics, this yields
\[
\dot x = \frac{1}{cx}[acx^2 + b(2a-c)x - b^2]\,.
\]
It is easy to see that under these dynamics, as $t\to -\infty$,
$x(t)\to\bar x$
with
\[
\bar x = \frac{b}{2ac}\left[-2a+c+\sqrt{4a^2+c^2}\right]\,.
\]

Notice also that the formula for $u$ in (\ref{singularm}) is increasing in $x$,
as is easily seen. Moreover, on the singular arc, $x\in[\hat x,\bar x]$, thus
on the one hand, $u \ge u^\sigma(\hat x)=0$, and on the other hand after
a straightforward calculation:
\[
u \le u^\sigma(\bar x) = \frac{1}{2}-\frac{\sqrt{4a^2+c^2}-2a}{2c} < \frac{1}{2}\,.
\]
Hence this is a feasible control.

The tributaries of this singular arc are investigated together with those of the
cooperative case in the appendix.

It is worth mentioning that, at the final point $(\hat t,\hat x)$ of this new
singular arc, we have for the singular control $u=0$. Hence the singular arc
has the same slope as the switch line, (and as the primary tangent to it). The
singular arc and the switch line form together a ``smooth'' curve. Such a
feature is known to be the rule in zero-sum two-person differential games.
The next section explains why.

\section{Further remarks}
\subsection{A zero-sum game formulation of Wardrop equilibria}

\begin{figure}[t]
$$\includegraphics[scale=1.]{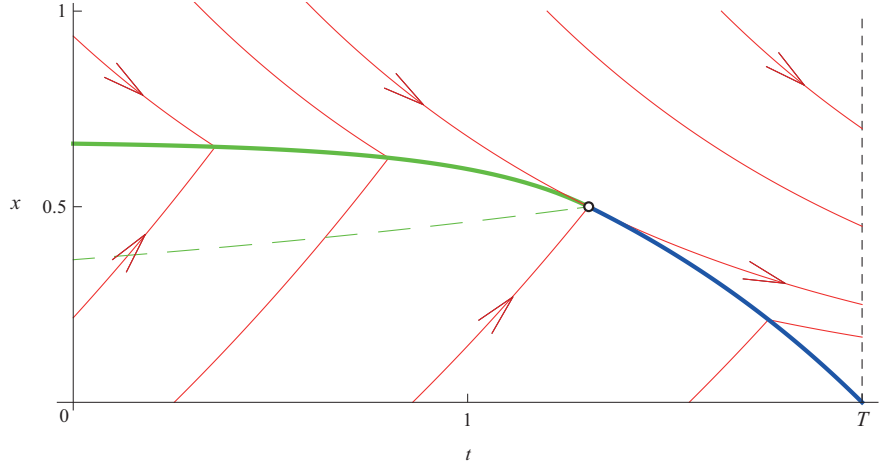}$$
\caption{The field of trajectories of the ESS. Green solid line: the singular arc.
Dotted green line: the singular arc of the cooperative solution.}
\label{f:1}
\end{figure}

Let a ``resident'' population ``chooses'' its behaviour, or strategy, $u$ in
some strategy
space $\mathcal{U}$, while a potential mutant ---or cheater, or free-rider---
would use a strategy $u_m\in\mathcal{U}$. Let $J(u,u_m)$ and $J_m(u,u_m)$
be the fitness of the resident and mutant population respectively under this
scenario. By hypothesis, if the ``mutant'' behaves as the resident, there
is no difference between them in terms of fitness, hence
\begin{equation}\label{nomut}
J(u,u) = J_m(u,u)\,.
\end{equation}
A strategy $u^\star\in\mathcal{U}$ is said to be \emph{uninvadable}, or
be a \emph{Wardrop equilibrium} \cite{war52}, if a Â mutant cannot do better
than a resident population using that strategy, i.e. if
\[
\forall u_m\in\mathcal{U}\,,\quad J_m(u^\star,u_m) \le J(u^\star,u_m)\,.
\]
Let
\begin{equation}\label{DeltaJ}
\Delta J(u,v) = J_m(u,v)-J(u,v)\,.
\end{equation}
Equivalently, the uninvadability property can be stated as
\[
\sup_{u_m\in\mathcal{U}}\Delta J(u^\star,u_m) \le 0\,.
\]
However, it follows from (\ref{nomut}) that
$\Delta J(u^\star,u^\star) = 0$, so that the supremum above
is always non-negative. Hence, a strategy $u^*$ is uninvadable, or a Wardrop
equilibrium, if and only if
\[
\max_{u_m\in\mathcal{U}}\Delta J(u^\star,u_m) = \Delta J(u^\star,u^\star)=0\,.
\]
And since this supremum is always non-negative, it follows that there
exists a Wardrop equilibrium if and only if
\[
\min_{u\in\mathcal{U}}\sup_{u_m\in\mathcal{U}}\Delta J(u,u_m) = 0\,.
\]
Therefore, the search for a Wardrop equilibrium can be performed as follows :
solve the zero-sum two-person game $\min_u\max_{u_m}\Delta J(u,u_m)$,
and if the min exists, check whether the min-sup reached yields a Value zero.
If yes, the minimizing $u$ is a Wardrop equilibrium. If furthermore the
maximum in $u_m$ is unique, then this is also an evolutionary stable strategy,
or ESS \cite{MSpri73,MS82}. This is the case in the above problem.

\subsection{Homogeneous problems}
The following piece of theory can undoubtly be deduced from Noether's
theorem \cite{noe18}. We offer a simple direct derivation. It accounts for
the transformation we used via the introduction of the variable $x$.

\subsubsection{The problem}
Let a control problem be given by a state $y\in\IR^n$, a control
$u\in \sU$, dynamics
\[
\dot y = f(y,u), \quad y(0)=y_0,
\]
with all regularity and growth assumptions to guarantee existence of
the solution for all measurable control functions
$u(\cdot) \in \mathcal{M}(\IR_+\to\sU)$.

A terminal condition $\cT(y)=0$, $\cT(\cdot)\in C^1(\IR^n)$, meaning that
$t_1 = \inf\{t\mid \cT(y(t))=0\}$, and for simplicity we assume that all
feasible trajectories are transverse to that differentiable manifold.
A criterion to be minimized by the control action is given, of the
form
\[
J(y_0,u(\cdot)) = K(y(t_1))+\int_0^{t_1}\!L(y(t),u(t))\,\d t\,.
\]

\paragraph{Assumptions} \mbox{}
\begin{enumerate}
\item The functions $f$, $\cT$, $K$ and $L$ are homogeneous of degree
one in $y$.
\item The state variable $y_n$ remains strictly positive for all
control functions.
\end{enumerate}
It follows from these assumptions that if $\phi$ stands for $f$, $\cT$,
$K$ or $L$, and if $y^{(n-1)}$ stands for the first $n-1$ components of
$y$, there exists a function of $n-1$ variables noted $\tilde\phi$
such that
\[
\phi(y) = y_n\tilde\phi\left(\frac{1}{y_n}y^{(n-1)}\right)
\]
Let also $x = (1/y_n)y^{(n-1)}$. Then, we claim
\paragraph{Theorem} \mbox{}
{\it
\begin{enumerate}
\item The reduced state $x(\cdot)$ obeys the dynamics
\[
\dot x = \tilde g(x,u) = \tilde f^{n-1}(x,u) - x\tilde f_n(x,u)\,,
\]
\item The Value function $V(y)$ of this problem is
homogeneous :
\[
V(y) = y_n\tV\left(\frac{1}{y_n}y^{(n-1)}\right)
\]
\item The function $\tV$ solves the Hamilton Jacobi equation
of the problem
\[
\tilde J(x_0) =
y_n(T)\widetilde K(y(T))+\int_0^{t_1}y_n(t)\tilde L(x(t),u(t))\,\d t\,.
\]
where
\[
y_n(t) := y_n(0)\exp\left(\int_0^t\tilde f_n(x(s),u(s))\,\d s\right).
\]
\end{enumerate}
}
\subsubsection{Proof of the theorem}
\paragraph{Dynamics}
Direct calculations result in
\[
\dot x_i = \frac{1}{y_n}f_i(y,u) - \frac{y_i}{y_n^2}f_n(y,u) =
\tilde g_i(x,u)\,.
\]
\paragraph{Value function}
It is completely elementary to see that if the data are globally
homogeneous of degree one in $x$, so is the Value function.
Therefore, there exists $\tV$ such that $V(x) = x_n\tV((1/x_n)x^{n-1}$.
\paragraph{Hamilton Jacobi equation}
Easy calculations yield
\[
\parder{V}{y_i}(y) = \parder{\tV}{x_i}(x), \> i=1,\ldots,n-1,
\]
and
\[
\parder{V}{y_n}(y) = \tV(x) - \sum_{i=1}^{n-1}x_i\parder{\tV}{x_i}(x).
\]
Substitute these in the standard Hamilton Jacobi equation
\[
\max_u\left[\parder{V}{y}f(y,u) + L(y,u)\right] = 0\,,
\]
to get (with $x:= (1/y_n)y^{(n-1)}$)
\[
\max_u\left[\parder{\tV}{x}(x) f^{n-1}(y,u) +
f_n(y,u)\left(\tV(x) - \parder{\tV}{x}(x)x\right) + L(y,u)\right]=0\,.
\]
If $y_n \neq 0$, we may write this as
\[
\max_uy_n\left[\tilde f_n(x,u)\tV(x) + \parder{\tV}{x}(x)\tilde g(x,u)
+ \tilde L(x,u)\right] = 0\,.
\]
If furthermore $y_n(t) > 0$, then we may divide through by $y_n$
without changing the $\max$ operator. We recognize the Hamilton
Jacobi Bellman equation of the problem with an exponential discount
factor $\exp\int\tilde f_n(x,u)$. However, multiplied by the positive
constant $y_n(0)$, this exponential is indeed $y_n(t)$.

The terminal condition $\cT(y)=0$ translates into
$\widetilde \cT(x) = 0$, and the boundary condition of the H.J.B.
equation into $\forall x : \widetilde\cT(x)=0$,
$\tV(x) = \widetilde K(x)$.

\subsection{Mixed strategy vs mixed population}
We want to investigate here the possible meanings of the mixed, singular
strategy. Recall that in our model, $u=0$ means reproducing, while $u=1$
means feeding. What is the meaning of an intermediary strategy $u\in(0,1)$ ?

One possible interpretation is that each individual in the population spends
some time eating and some time reproducing in a fast cycle, as compared to
the horizon of the problem ---here a year--- and to the characteristic frequency
of the dynamics, here $a$ and $c$. This thus a monomorphic population
agreeing on a mixed strategy for all individuals.

In linear problems, it is possible to think of an intermediary $\tilde u$ as
representing a polymorphic population, where a fraction $\tilde u$ uses the
strategy $u=1$ while a fraction $1-\tilde u$ uses the strategy $u=0$.
However, our criterion here is nonlinear. Therefore further investigation is
needed.

Let the population be represented by a measured space $\Om$, with a
positive measure $\mu$ of total mass $\mu(\Om)=1$. (It is thus a probability
measure.) Each individual is represented as an $\om\in\Om$. Assume also
that each individual $\om$ has to pick, at each instant of time, a control
$v(\om,t)\in\{0,1\}$. We need to assume that $\om\mapsto v(\om,t)$ is
$\mu$-measurable for (allmost) all $t\in [0,T]$, and that, for (allmost) all
$\om\in\Om$, $t\mapsto v(\om,t)$ is piecewise constant, equal to $0$ or
$1$ on time intervals (hence Lebesgue measurable).

Each individual has an energy $\pi(\om,t)$ at time $t$. Its dynamics are
\[
\dot\pi(\om,t) = -a\pi(\om,t) + b\n(t)v(\om,t)\,.
\]
The rsource depletion rate is
\[
\dot\n = -c\n\int_\Om\!v(\om,t)\,\d\mu(\om)\,.
\]
So, if we set
\[
u(t) = \int_\Om\!v(\om,t)\,\d\mu(\om) = \mu(\{\om\mid v(\om,t)=1\})
\]
and
\[
p(t) = \int_\Om\!\pi(\om,t)\,\d\mu(\om)\,,
\]
we find equations (\ref{pdot}) and (\ref{ndot}), as it should.

The number of offspring produced by an individual is
\[
f(\om) = \int_0^T\!(1-v(\om,t))\pi(\om,t)\,\d t\,,
\]
and for the whole population (using Fubini's theorem)
\[
F = \int_\Om\!f(\om)\,\d\mu(\om) =
\int_0^T\!\left[\int_\Om(1-v(\om,t))\pi(\om,t)\,\d\mu(\om)\right]\d t
\]
which would coincide with the formula (\ref{J}) for $J$ only if $\pi(\om,t)$
and $v(\om,t)$ were probabilistically independant, an impossible situation
if each player uses a constant control on nonzero, measurable, time intervals.
We would approach such an independance if we were to assume that the time
intervals during which the individuals'controls are constant are extremely
short, and after each, the set of individuals using $u=1$, say, is drawn at
random, with probability $u(t)$, independantly of $\pi(\om,t)$. But this is
reconstructing a monomorphic population of mixed players.

\appendix
\section{Switch function on feeding tributaries}
\subsection{All tributaries with $u=1$}
Let $t_s$, $x_s$, $\lambda_s$, $\mu_s$ stand for the variables $t$, $x$,
$\lambda$  and $\mu$ at the point where a trajectory with $u=1$ meets
either the switch line or a singular arc, either cooperative or uninvadable.
To simplify the notations, let also
\[
\alpha := \e^{-a(t_s-t)}\,,\qquad \gamma:= \e^{-c(t_s-t)}\,.
\]
On this trajectory, we have
\begin{eqnarray*}
x(t) &= &\left(x_s+\frac{b}{c-a}\right)\alpha-\frac{b}{c-a}\,,\\
\lambda(t) &= &\lambda_s\alpha\,,\\
\mu(t) &= &\mu_s\gamma + \frac{b}{c-a}(\alpha-\gamma)
\end{eqnarray*}
We want to investigate the sign of
\[
\sigma = b\lambda - c\mu - x\,,
\]
remembering that in all cases to be investigated, $\sigma(t_s)=0$.
We have already noticed that
\[
\dot\sigma = c\sigma + ab\lambda + ax - b\,.
\]
Hence
\begin{equation}\label{sigma}
\sigma = -\int_t^{t_s}\e^{-c(\tau-t)}(ab\lambda(\tau)+ax(\tau)-b)\,\d\tau\,.
\end{equation}
We wish to show that $L(\tau) := ab\lambda(\tau)+ax(\tau)-b < 0$, thus
proving that $\sigma$ is positive from $t=0$ to $t_s$. Substituting the
explicit values of $x$, $\lambda$ and $\mu$, we get
\[
L(t) = a\left[b\lambda_s\alpha + x_s\frac{\gamma}{\alpha}
+ \frac{b}{c-a}\left(\frac{\gamma}{\alpha}-\frac{c}{a}\right)  \right]
\]

\subsection{Tributaries of the switch line}
On the switch line, $x_s=b\lambda_s$ and $\mu_s=0$. Hence
\[
L = ab\left[\lambda_s\left(\alpha+\frac{\gamma}{\alpha}\right)
+\frac{1}{c-a}\left(\frac{\gamma}{\alpha}-\frac{c}{a}\right)
\right]
\]
We notice that $L(t_s)=ab(2\lambda_s - 1/a)$. But we know that on the switch
line $\lambda_s \le 1/2a$. Hence $L(t_s)\le 0$. Let us investigate its time
derivative:
\[
\frac{1}{ab}\dot L(t) =
\lambda_s\left(a\alpha+(c-a)\frac{\gamma}{\alpha}\right) +
\frac{\gamma}{\alpha}=
a\lambda_s\alpha +\left((c-a)\lambda_s+1\right)\frac{\gamma}{\alpha}\,.
\]
Hence $\frac{1}{ab}\dot L > (-a\lambda_s + 1)\gamma/\alpha > 0$, and
therefore, $L(t) \le 0$ for allt $t \le t_s$.

\subsection{Tributaries of the singular arcs}
We now have
\[
L(t_s) = a(b\lambda_s+x_s)-b\,.
\]
We claim that $x_s < b/a$. As a matter of fact, on the singular arc of the
cooperative solution, $x_s \le b/2a$. On the singular arc of the uninvadable
solution,
\[
x_s \le \bar x = \frac{b}{2a}\left[1 + \frac{\sqrt{4a^2+c^2}-2a}{c} \right]\,.
\]
It is easy to chack that $\sqrt{4a^2+c^2}-2a < c$, so that indeed, $x_s<b/a$,
and therefore, $L(t_s) < 0$. Again, compute its time derivative:
\[
\dot L(t) = a\left[ ab\lambda_s\alpha +
\left((c-a)x_s+b\right)\frac{\gamma}{\alpha}\right] > 0
\]
again thanks to $b-ax_s > 0$. So, again, we may conclude that for all
$t\le t_s$, $L(t) < 0$.

In all three cases, we conclude with the help of formula (\ref{sigma}) that on
these trajectories, for all $t\le t_s$, $\sigma(t) > 0$.

\subsection*{Acknowledgement}

This work is supported by Agropolis Foundation and the R\'eseau National des Syst\`emes Complexes (RNSC) under the ModPEA Project.

\addcontentsline{toc}{section}{\textit{References}}
\bibliographystyle{unsrt}

\end{document}